# Fractal structure of iterative time profiles

Marek Berezowski, Katarzyna Bizon
Silesian University of Technology, Institute of Mathematics, Gliwice, Poland
E-mail: marek.berezowski@polsl.pl

## Abstract

The scope of the paper is the theoretical analysis of the time rate in which a dynamical system reaches a stable stationary state or stable oscillations. The method used for the analysis is based on the so-called iterative time profiles, demonstrating a chaotic and fractal nature of some of the profiles. The results were presented in the form of two and three-dimensional graphs.

## 1. Introduction

The dynamics of a reactor was discussed, for example, in [1], [2], [3], [4] and [5], where it was demonstrated that the variables may be constant in the steady state or may oscillate in the form of rectangular time wave. The oscillations may be more or less complex, depending on the periodicity of the time series of the above-mentioned variables of a given state. Due to an asymptotic nature of the changes occurring in the reactor, the time rate required for the system to reach a steady state is, from the theoretical point of view, infinitely long. In practice, however, it may be assumed that the reactor is in the steady state when the values of its variables are relatively close to the values of the theoretical state. Thus, the time rate required for the system to reach the steady state is limited and depends on the values of the reactor parameters as well as on the values of temperature and concentration of the fluid flux at the initial moment. The reactor subjected to the theoretical analysis with external recycle [1], [3] and [5]. The time required for the reactor to reach the steady state was illustrated by means of various types of two and three-dimensional iterative time profiles [6]. The profiles were demonstrated to have a chaotic nature and fractal structure.

## 2. The model and its iterative time profiles

The differential model of a reactor was analyzed in the following forms ([1], [3]):

$$\frac{d\alpha_{k+1}(\xi)}{d\xi} = (1-f)\phi[\alpha_{k+1}(\xi), \Theta_{k+1}(\xi)] \tag{1}$$

$$\frac{d\Theta_{k+1}(\xi)}{d\xi} = (1-f)\{\phi[\alpha_{k+1}(\xi), \Theta_{k+1}(\xi)] + \delta[\Theta_H - \Theta_{k+1}(\xi)]\} \tag{2}$$

$$\alpha_{k+1}(0) = f\alpha_k(1); \quad \Theta_{k+1}(0) = f\Theta_k(1) \tag{3}$$

$$\phi(\alpha, \Theta) = Da(1-\alpha)^n \exp\left(\gamma\frac{\beta\Theta}{1+\beta\Theta}\right). \tag{4}$$

The purpose of the analysis was to determine the impact of the initial values of $\alpha_0(1)$ and $\Theta_0(1)$ on the number of discrete time steps $N$ required for the reactor to reach a stable steady state, that is stable stationary point or stable periodic orbit [6]. It was assumed that the steady state for the reactor is the state when the distance between the trajectory generated by the model (1), (2), (3) and (4) and the fixed point $\alpha_s(1)$, $\Theta_s(1)$ is not longer than the one determined by the following condition:

$$\left(\left|\frac{\alpha(1) - \alpha_s(1)}{\alpha_s(1)}\right| + \left|\frac{\Theta(1) - \Theta_s(1)}{\Theta_s(1)}\right|\right) * 100\% < \varepsilon. \tag{5}$$

The analysis was based on the catastrophic set shown in Fig. 1 [1] and [3], generated for the following values of the parameters: $Da = 0.15$, $n = 1.5$, $\gamma = 15$, $\beta = 2$, $\delta = 3$. Crossing of the lines SN, FB or HB of this set leads to the change in the multiplication factor (saddle-node bifurcation), generation of jump oscillations (flip bifurcation) or generation of quasi-periodic oscillations (Hopf bifurcation), respectively. Thus, right to the zone marked by the FB line, the temperature and concentration of the fluid flux are constant in the steady state, practically for any value of the recycle coefficient within the range of $0 < f < 1$. Accordingly, what we are dealing with are fixed stationary points. But, within the boundaries of the zone marked by the FB line, the temperature and concentration oscillate in the steady state in a jump way [3]. Accordingly, we are dealing with fixed points of oscillation solutions.

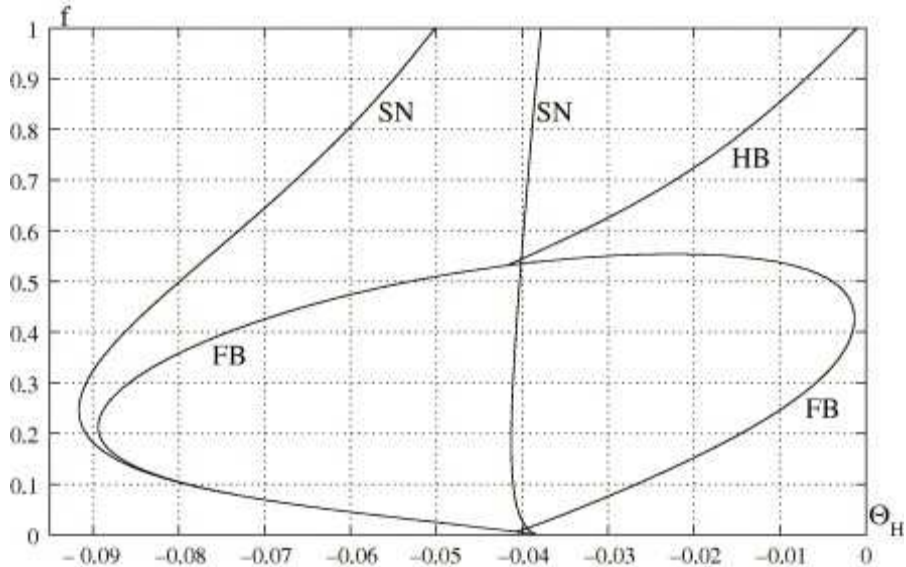

Fig. 1. Catastrophic set.

By analyzing the eigenvalues $\alpha$ of the model (1), (2), (3), (4) and (5), it is possible to determine the "*power*" by means of which a given fixed point attracts trajectories. The most powerful attractor is the point for which all the eigenvalue modules are equal to zero, whereas the least powerful point is the one for which the module of at least one eigenvalue is equal to 1. The eigenvalues of $k$-periodic solution may be derived from the following relation:

$$\bar{y}_k = f^k \prod_{j=1}^{k} e^{\int_0^1 \bar{\bar{J}}_j d\xi} \bar{y}_0 = f^k \prod_{j=1}^{k} \overline{\overline{M}}_j \bar{y}_0 = f^k \overline{\overline{A}}_k \bar{y}_0 \tag{6}$$

where $\bar{\bar{J}}_j$ is Jacobi matrix of the right sides of Eqs. (1) and (2), whereas $\bar{y}$ is the state vector of the linear approximation of the equations. The monodromy matrix $\bar{\bar{M}}_j$ may be derived from the following equation:

$$\frac{d\bar{\bar{M}}_j}{d\xi} = \bar{\bar{J}}_j \bar{\bar{M}}_j \qquad (7)$$

assuming, however, that $\bar{\bar{M}}_j(\xi = 0)$ is an identity matrix. The searched eigenvalues of the model (1), (2), (3), (4) and (5) are the eigenvalues of matrix $f^k \bar{\bar{A}}_k$.

Furthermore, assuming that $\Theta_H = -0.001$ (right to the zone marked by FB), the bifurcation graph shown in Fig. 2 was designated, from which it may be indicated that the model of the reactor has a stationary solution within the full variation range of $0 < f < 1$. Assuming that $k = 1$, the dependence between the modules of eigenvalues $\lambda_i$ ($i = 1, 2$) and the recycle coefficient was shown in Fig. 3. As seen in the figure, the maximal eigenvalue occurs for $f = 0.427$. Accordingly, it is for this value of the recycle coefficient that the stationary state of the reactor is a weak attractor of the trajectories. This means a long way for the trajectory to lead the reactor to the state of dynamical equilibrium. In such case the number of recurrence steps $N$ is big. To substantiate the above deliberation, the tree of the recurrence solutions is shown in Fig. 4, where, for a given value of $f$, different values of $N$ were marked, corresponding to different initial values of the conversion degree $0 < \alpha_0(1) < 1$. The maximal value corresponds to the maximal eigenvalue from Fig. 3. The calculations were based on the assumption that $\Theta_0(1) = 0.2$.

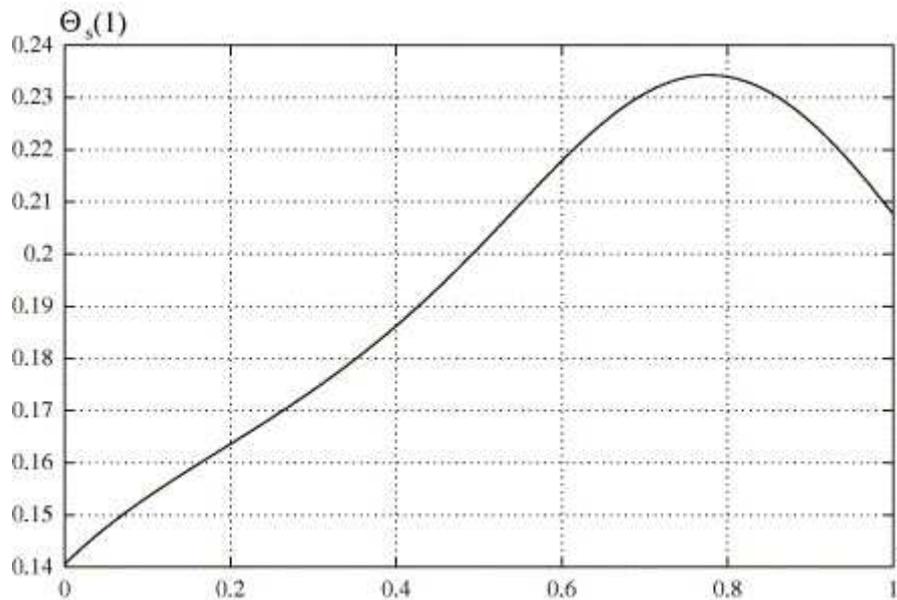

Fig. 2. Bifurcation diagram, $\Theta_H = -0.001$.

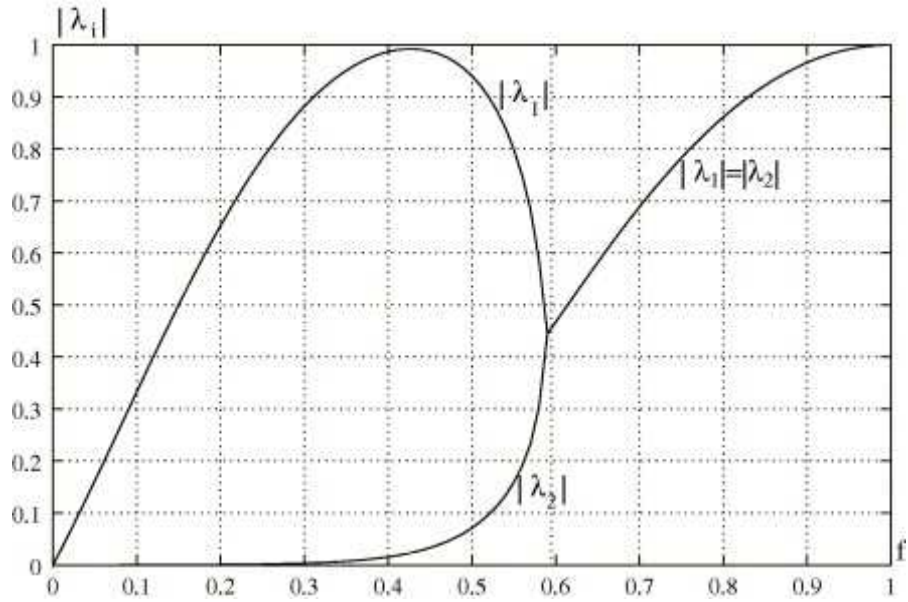
Fig. 3. The eigenvalues of the model, $\Theta_H = -0.001$.

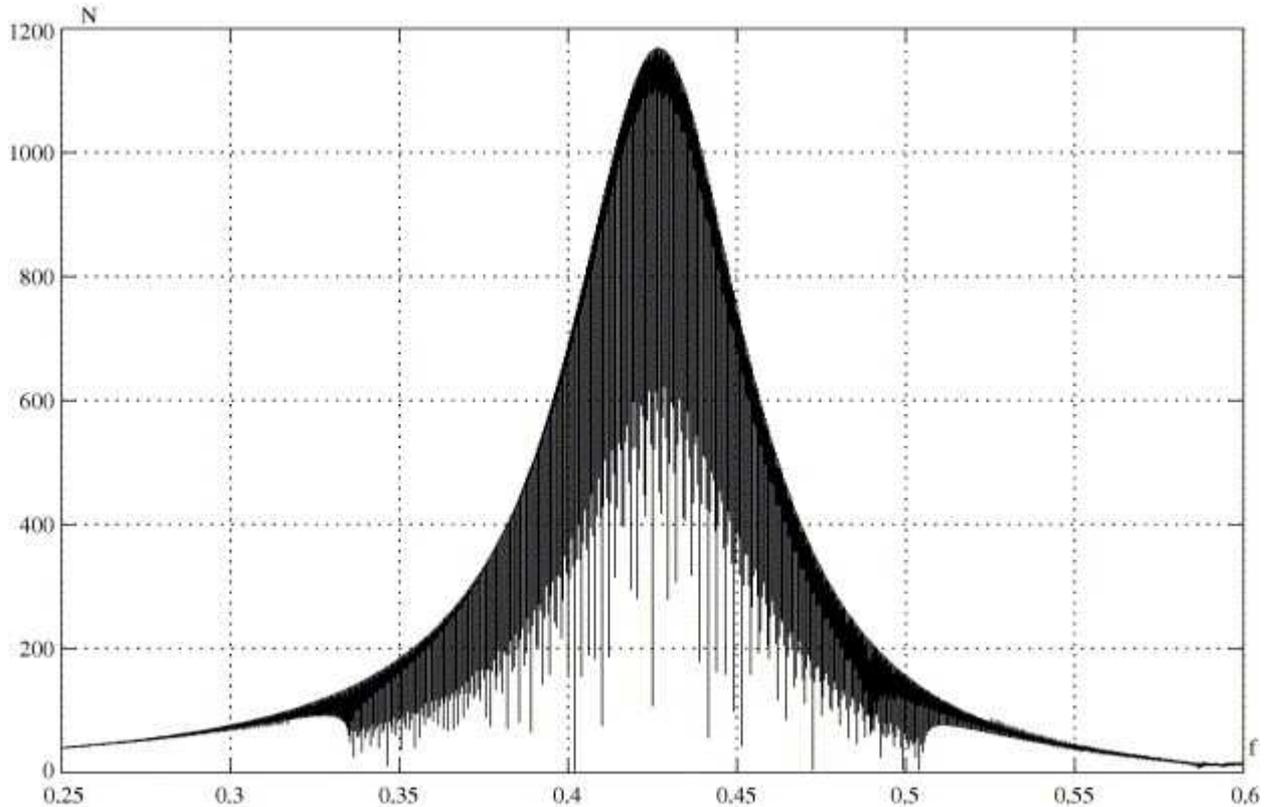
Fig. 4. Tree of iterative time profiles, $\Theta_H = -0.001$.

Assuming that $\Theta_H = -0.002$ (the internal zone marked by the line FB), the bifurcation diagram shown in Fig. 5 was designated, where, apart from stationary solutions, twoperiodic oscillation solutions in the range of $0.3795 < f < 0.4695$ may be observed. Assuming that $k = 2$, the dependence between the modules of eigenvalues $\lambda_i$ ($i = 1, 2$) and the recycle coefficient was shown in Fig. 6. As seen in the figure, the maximal value from Fig. 3 corresponds (with big accuracy) to the minimum from Fig. 6. Because the minimum lies in the above-mentioned variation range of parameter $f$, it concerns an oscillation solution. For the value of the recycle coefficient $f = 0.427$ the twoperiodic orbit is a strong attractor of the trajectories. This means a short way for the trajectory to lead the reactor to the stable temperature and concentration oscillations. In such case the number of recurrence steps $N$ is small. At the ends of the variation range of parameter $f$, that is, for $f = 0.3795$ and $f = 0.4695$, the eigenvalue is equal to 1. Accordingly, at this particular position the fixed point is the least powerful attractor of the trajectories. To substantiate the above deliberation, another tree of recurrent solutions is shown in Fig. 7, where, the minimum corresponds to the minimal eigenvalue from Fig. 6 in the range of $0.3795 < f < 0.4695$. The two maximums that are clearly seen in the figure designate the range of oscillations FB, equal to the range presented in Fig. 5 and Fig. 6.

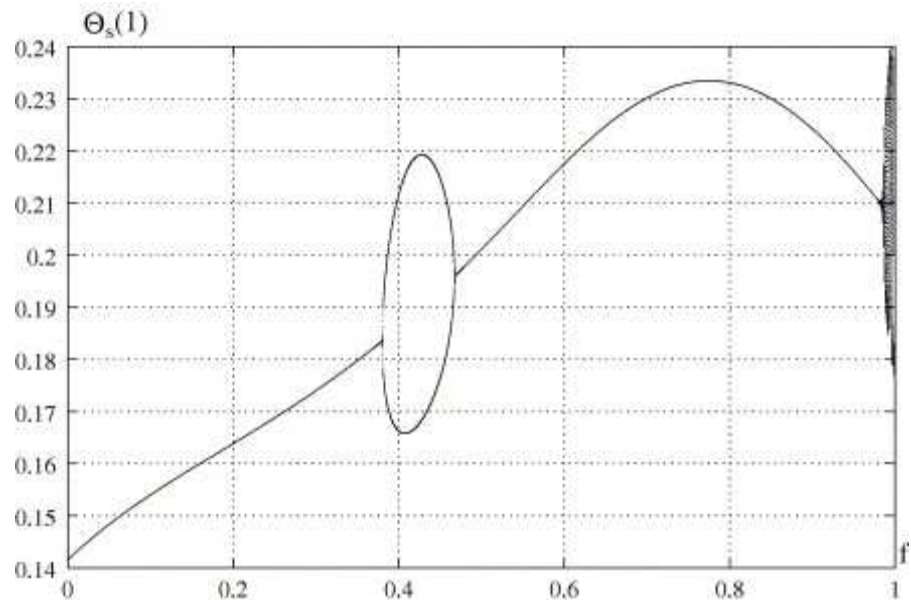

Fig. 5. Bifurcation diagram, $\Theta_H = -0.002$.

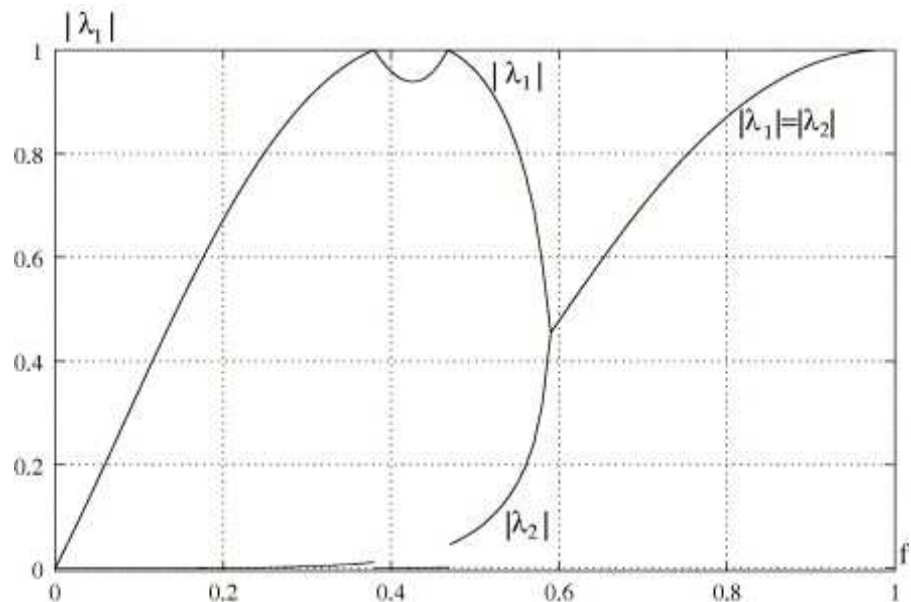

Fig. 6. The eigenvalues of the model, $\Theta_H = -0.002$.

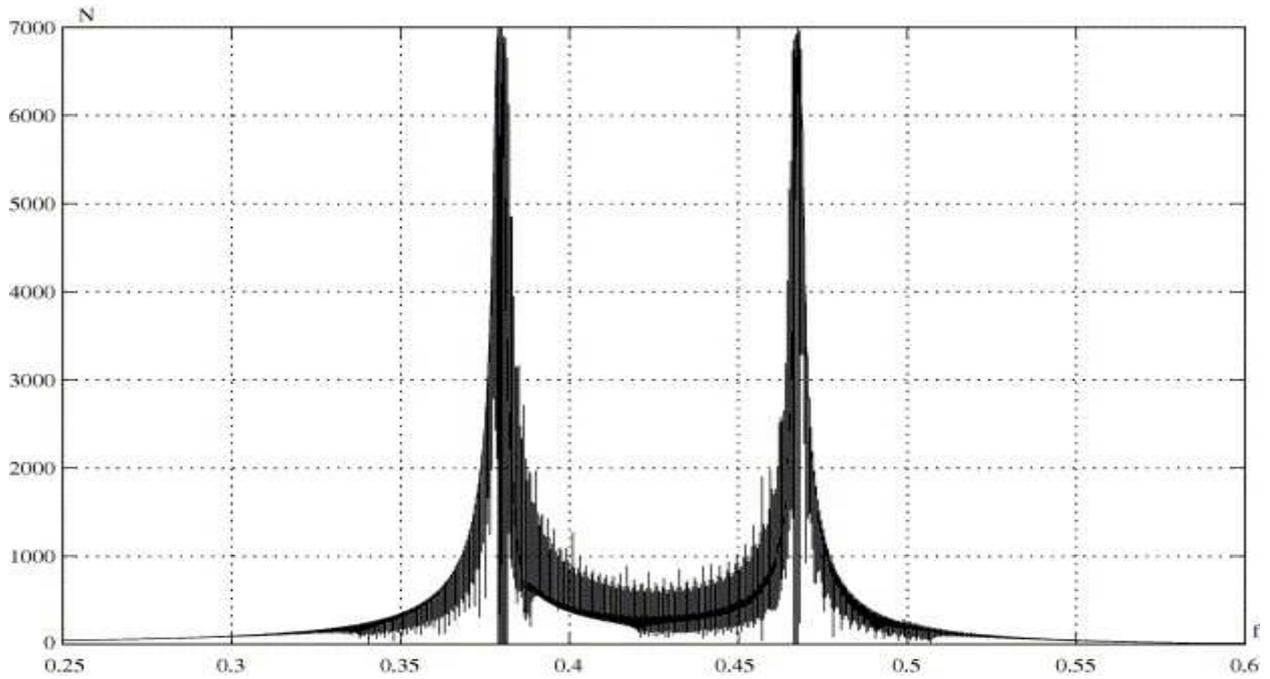

Fig. 7. Tree of iterative time profiles, $\Theta_H = -0.002$.

By decreasing the value of the dimensionless temperature of the cooling medium the model generates another FB bifurcation, which leads to the doubling of the period of the present oscillation solution. In consequence, for $\Theta_H = -0.012$, fourperiodic oscillations in the range of $0.381 < f < 0.430$ are generated. To substantiate the above deliberation, another tree of recurrent solutions is shown in Fig. 8, where, the three minimums correspond to strong attraction by a fourperiodic attractor, and the four maximums to FB bifurcation points—which are weak attractors for the trajectories. Decreasing further the value of $\Theta_H$ other FB bifurcations are obtained, generating the successive doubling of the period [3]. Each of the bifurcations offers $(k-1)$ minimums, i.e. $(k-1)$ strong attractors of the trajectories and $k$ maximums, i.e. $k$ weak attractors.

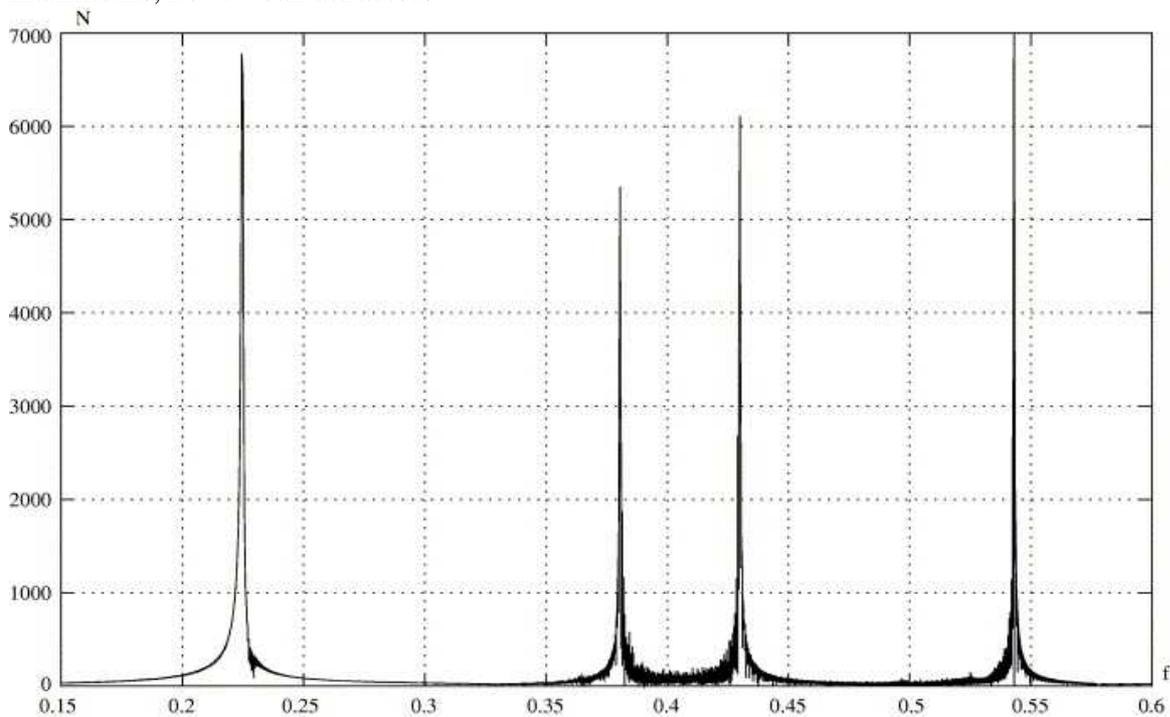

Fig. 8. Tree of iterative time profiles, $\Theta_H = -0.012$.

For example, assuming that $\Theta_H = -0.001$, $f = 0.427$ and $\Theta_0(1) = 0.2$ an iterative time profile demonstrating the impact of the initial value of the conversion degree on the number of recurrence steps $N$ required for the reactor to reach a stable stationary state is shown in Fig. 9a with the accuracy of $\varepsilon = 0.001\%$. For $\alpha_0(1) = 0.559$ the visible peak directed downwards certifies that the initial condition of the reactor is located closest to the stationary attractor.

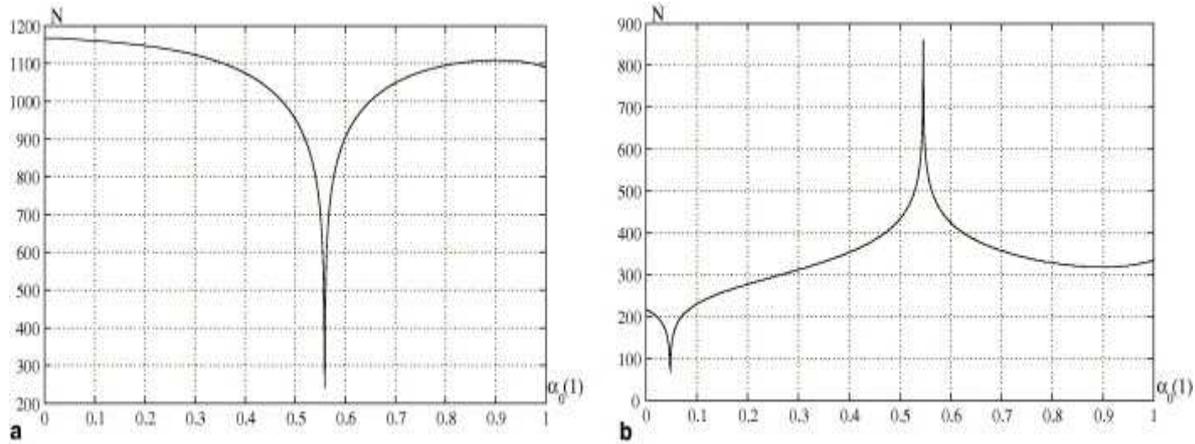

Fig. 9. (a) Iterative time profile, $\Theta_H = -0.001$, $f = 0.427$, $\Theta_0(1) = 0.2$. (b) Iterative time profile, $\Theta_H = -0.002$, $f = 0.427$, $\Theta_0(1) = 0.2$.

Likewise, in Fig. 9b the iterative time profile for $\Theta_H = -0.002$, $f = 0.427$ and $\Theta_0(1) = 0.2$ demonstrates the impact of the initial value of the conversion degree $\alpha_0(1)$ on the number of recurrence steps $N$ required for the reactor to reach a stable limit cycle. For $\alpha_0(1) = 0.546$ the visible peak directed upwards certifies that the initial condition of the reactor is located farthest from the twoperiodic attractor. In view of an insignificant discrepancy in the values of $\Theta_H$ in both of the above-mentioned cases, the peaks are also insignificantly displaced from each other, yet directed oppositely. The change in the direction of the peak is a result of crossing the bifurcation line FB on the graph from Fig. 1. This regularity also concerns the next bifurcations, leading to the successive doubling of the period of the oscillation solutions. In Fig. 9b the visible peak directed downwards [$\alpha_0(1) = 0.048$] certifies that the initial condition of the reactor is located closest to the twoperiodic attractor. Further decrease in the value of parameter $\Theta_H$ leads to the occurrence of successive FB bifurcations, which, in consequence, leads to the generation of more and more peaks on the iterative time graph. The number of the peaks may increase very rapidly. In Fig. 10 the iterative time profile for $\Theta_H = -0.012$, $f = 0.427$ and $\Theta_0(1) = 0.2$ has a chaotic nature in the range of $0.4 < \alpha_0(1) < 0.42$. This chaos is revealed by a great sensitivity of this part of the graph to the infinitesimal change in the position of a fourperiodic attractor, as well as by a disordered cloud of the points on Poincaré cross-section (Fig. 11). The number of iteration steps $N$ corresponding to the $j$th and $(j + 1)$th peak from Fig. 10 is marked on the coordinates of the cross-section. Moreover, the graph in Fig. 10 was found to have a fractal nature, which was demonstrated on three-dimensional iterative graphs ( Fig. 12 and Fig. 13).

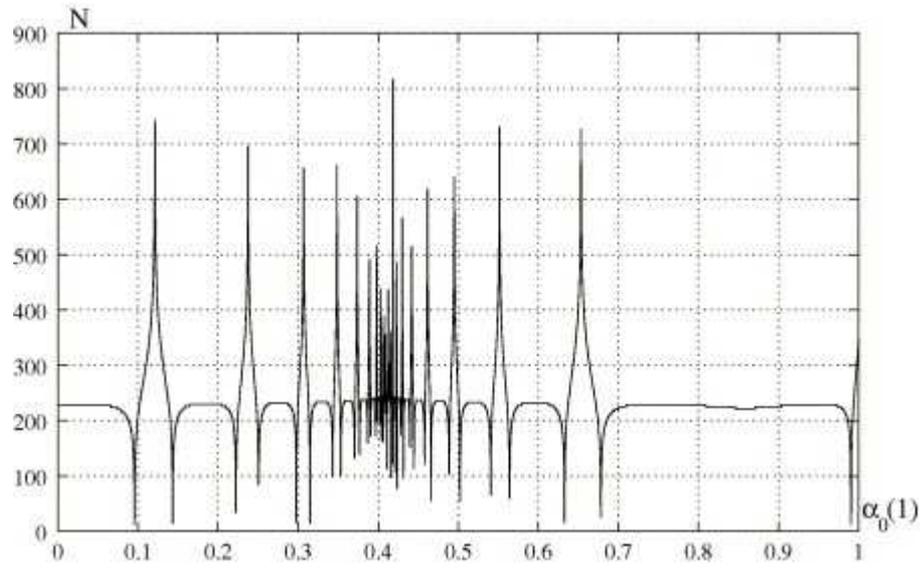
Fig. 10. Iterative time profile, $\Theta_H = -0.012, f = 0.427, \Theta_0(1) = 0.2$.

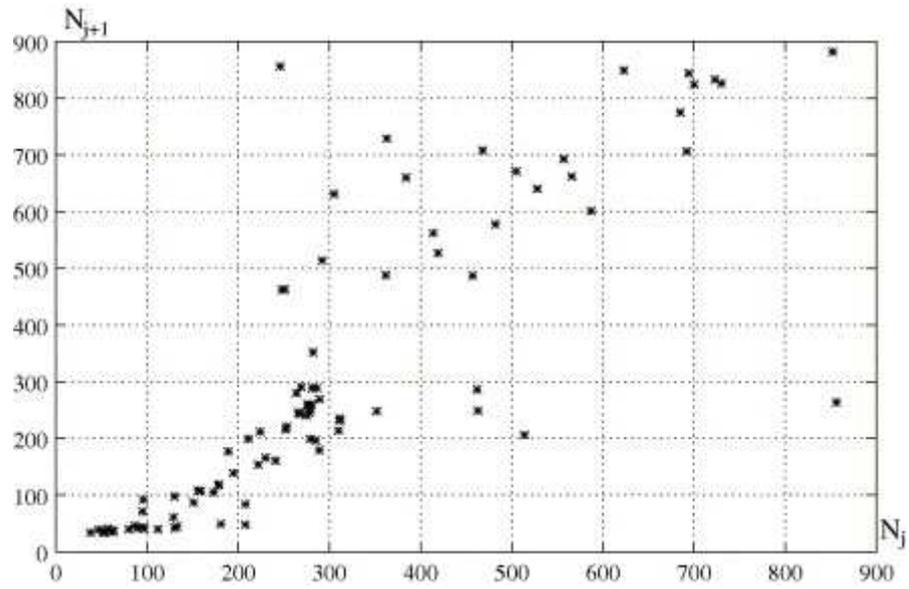
Fig. 11. Poincaré cross-section, $\Theta_H = -0.012, f = 0.427, 0.4 < \alpha_0(1) < 0.42$.

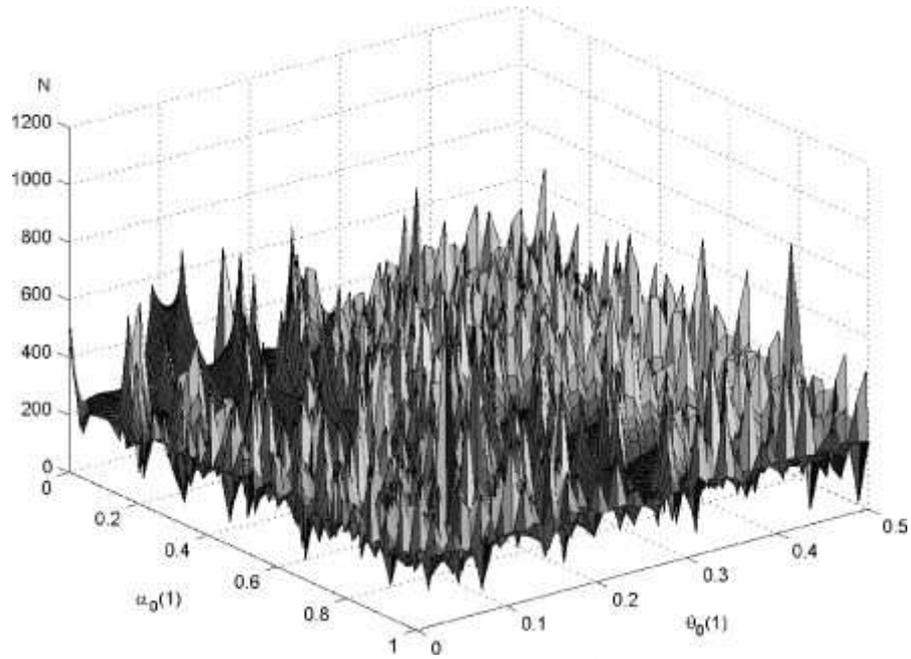

Fig. 12. Three-dimensional iterative time profile, $\Theta_H = -0.012, f = 0.427$.

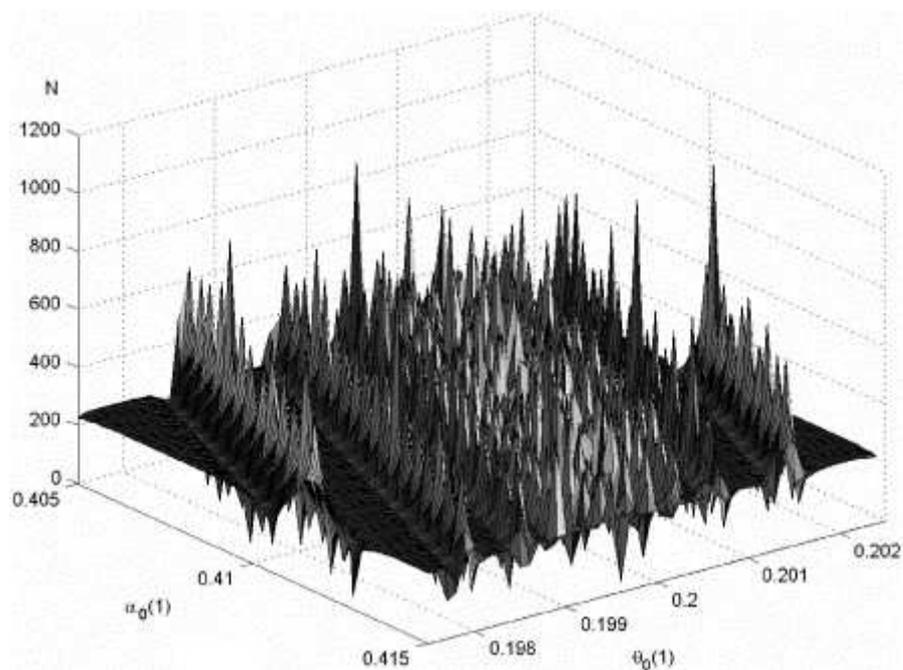

Fig. 13. Fragment of the profiles from Fig. 12.

Other exemplary fractal diagrams of iterative time profiles of the discussed reactor are to be found in [7].

## 3. Concluding remarks

The scope of the paper was a theoretical analysis of the dynamics of a reactor. Particular attention was given to the impact of the process parameters and the initial values of the concentration and temperature of the fluid flux on the time rate required for the reactor to reach steady states. In the course of the research the so-called iterative time profiles were designated and their chaotic and fractal nature was shown. The results are useful from the view point of process technology, and interesting from the theoretical point of view. Due to substantial and significant mass effects and thermal effects caused by the chemical reaction, the discussed method makes it possible to optimize the design of the reactor system, whereas the results give an insight into the dynamics of reactors.